\documentclass[a4paper, 12pt]{article}

\usepackage[koi8-r,cp1251]{inputenc}
\usepackage[english,russian]{babel}
\usepackage{amsfonts,amssymb,amsmath, hyperref}
\usepackage[final]{epsfig}
\usepackage{graphicx}

\begin{document}

\begin{flushleft}
 УДК 517.983.23
\end{flushleft}

\begin{center}
{\bf О МНОГОМЕРНОМ ФУНКЦИОНАЛЬНОМ ИСЧИСЛЕНИИ БОХНЕРА-ФИЛЛИПСА\\
 \copyright 2009 г. А. Р. Миротин}\\
{\it Гомельский государственный университет им. Ф. Скорины,
Гомель}
\end{center}

\hspace{5mm}
\begin{center}
{\bf ON MULTIDIMENSIONAL BOCHNER-PHILLIPS FUNCTIONAL CALCULUS\\
 A. R. Mirotin}\\
{\it F. Scorina Gomel State University, Gomel}\\
amirotin@yandex.ru
\end{center}

\hspace{5mm}
 Развивается принадлежащее автору многомерное
функциональное исчисление генераторов полугрупп, основанное  на
классе функций Бернштейна нескольких переменных. Уточняется
основная теорема этого исчисления,  дается условие голоморфности
полугрупп, порождаемых операторами, возникающими в исчислении, и
для одномерного случая доказывается неравенство моментов для этих
операторов.

\hspace{5mm}

{\it Ключевые слова: полугруппа операторов, генератор
полугруппы, функциональное исчисление, функция Бернштейна}

\hspace{5mm}

The functional calculus of semigroup generators, based on the
class of Bernstein functions in several variables is developed,
the condition for holomorphy  of semigroups, generated by
operators which arisen in the calculus  is given, and in the
one-dimensional case the moment inequality for such operators is
proved.

\hspace{5mm}


\

The paper was published in [14].

\hspace{5mm}

{\it  1 Введение и предварительные сведения}

 Одномерное функциональное исчисление Бохнера-Филлипса
является существенной частью теории полугрупп операторов (см.,
например, [1 -- 3]) и находит важные применения в теории случайных
процессов. Основы многомерного исчисления были заложены автором в
[4 -- 6]. Целью настоящей работы является дальнейшее развитие
этого исчисления. Прежде всего, мы уточняем и дополняем основную
теорему из [6]. Далее,  указываются условия, при которых
операторы, возникающие в рассматриваемом исчислении, порождают
голоморфные полугруппы. Как следствие получается положительный
ответ на один вопрос Кишимото и Робинсона [1] для равномерно
выпуклых пространств. Наконец, для одномерного случая мы
доказываем полезное неравенство моментов, обобщающее неравенство
из [7].

Для формулировки результатов напомним необходимые понятия и факты
из [6].

{\bf 1. Определение [4]}. Будем говорить, что неположительная
функция $\psi\in C^\infty((-\infty;0)^n)$ {\it принадлежит классу
${\cal T}_n$}, если все ее частные производные первого порядка
абсолютно монотонны (функция из $C^{\infty}((-\infty;0)^n)$
называется {\it абсолютно монотонной,} если она неотрицательна
вместе со своими частными производными всех порядков).

Последнее условие на $\psi$ равносильно тому, что
$\partial^\alpha\psi\geq 0$ для любого мультииндекса $\alpha\ne
0$.

Очевидно также, что ${\cal T}_n$ есть конус относительно
поточечного сложения функций и умножения на скаляр.

Известно [6], что каждая функция $\psi\in {\cal T}_n$ допускает
интегральное представление ($s\in(-\infty;0)^n$)

$$\psi(s)=c_0+c_1\cdot s+\int\limits_{\Bbb{R}_+^n\setminus \{0\}}
(e^{s\cdot u}-1)d\mu (u),     \eqno(1.1)
$$

\noindent где $c_0=\psi(-0):=\lim\limits_{s\to -0}\psi(s)$, а
$c_1$ из $\Bbb{R}_+^n$ и положительная мера $\mu$ на
$\Bbb{R}_+^n\setminus \{0\}$ определяются однозначно (здесь и ниже
точкой мы обозначаем скалярное произведение в $\Bbb{R}_+^n$;
запись $s\to -0$ означает, что $s_1\to -0, \ldots, s_n\to -0$).

Всюду далее через $T_1, \ldots ,T_n$ будут обозначаться попарно
коммутирующие однопараметрические $C_0$-полугруппы (т. е. сильно
непрерывные на $\Bbb{R}_+$ полугруппы) в комплексном банаховом
пространстве $X$, удовлетворяющие условию $||T_j(t)||\leq
M\quad(t\geq 0;j=1, \ldots ,n; M={\rm const}\geq 1)$. Через $A_j$
обозначим генератор полугруппы $T_j$ с областью определения
$D(A_j)$ и положим $A=(A_1,\ldots ,A_n)$. Далее коммутирование
операторов $A_1,\ldots, A_n$   означает коммутирование
соответствующих полугрупп. Через ${\rm Gen}(X)$ мы будем
обозначать множество всех генераторов равномерно ограниченных
$C_0$-полугрупп в $X$, а через $I$ -- единичный оператор в $X$.

Операторнозначная функция $T(u):=T_1(u_1)\ldots
T_n(u_n))\quad(u\in \Bbb{R}_+^n)$ является
$C_0$-$n$-параметрической полугруппой, а потому линейное
многообразие $D(A):=\cap_{j=1}^nD(A_j)$ плотно в $X$ ([8], c. 98
-- 99).

{\bf 2. Определение [6]}. Определим значение функции $\psi$ из
${\cal T}_n$ вида (1.1) на наборе $A=(A_1,\ldots ,A_n)$ при $x\in
D(A)$ формулой
$$
\psi(A)x=c_0x+c_1\cdot Ax+\int\limits_{\Bbb{R}_+^n\setminus \{0\}}
(T(u)-I)xd\mu(u),     \eqno(1.2)
$$
где $c_1\cdot Ax:=\sum_{j=1}^nc_1^jA_jx. $

Пусть $\psi\in{\cal T}_n, t\geq 0$. Тогда функция
$g_t(s):=e^{t\psi(s)}$ будет абсолютно монотонной на
$(-\infty;0)^n$. Очевидно также, что $g_t(s)\leq 1.$ В силу
многомерного варианта теоремы Бернштейна-Уиддера (см., например,
[9], c. 281) существует такая единственная ограниченная
положительная мера $\nu_t$ на $\Bbb{R}_+^n,$ что при $s\in
(-\infty;0)^n$
$$
g_t(s)=\int\limits_{\Bbb{R}_+^n} e^{s\cdot u}d\nu_t(u)=({\cal
L}\nu_t)(-s),
$$
\noindent где ${\cal L}$ обозначает $n$-мерное преобразование
Лапласа.

{\bf 3. Определение.} Используя обозначения, введенные выше,
положим ($x\in X$)
$$
g_t(A)x=\int\limits_{\Bbb{R}_+^n} T(u)xd\nu_t(u)\eqno(1.3)
$$
(интеграл понимается в смысле Бохнера).

Очевидно, что $||g_t(A)||\leq M^ne^{t\psi(-0)}\leq M^n.$ Поскольку
$g_{t+r}(s)=g_t(s)g_r(s),$ то $\nu_t$ образуют сверточную
полугруппу ограниченных мер на $\Bbb{R}_+^n$. Поэтому
$g(A):t\mapsto g_t(A)$ есть равномерно ограниченная полугруппа
операторов на $X$. В частности, $g(A)$ есть $C_0$-полугруппа.

Введенные выше обозначения и ограничения далее будут применяться
без дополнительных пояснений.

\vspace{5mm}
 {\it 2 Основные результаты}

{\bf 4. Теорема}. {\it Замыкание оператора $\psi(A)$ существует и
является генератором  полугруппы $g(A)$ класса $C_0$, определенной
формулой \rm{(1.3)}.}

Доказательство. В [6] было доказано, что оператор  $\psi(A)$
замыкаем, и его расширением является генератор $G$
$C_0$-полугруппы $g(A)$.  Поэтому замыкание
$\overline{\psi(A)}\subseteq G$. Покажем, что здесь имеет место
равенство. Поскольку операторы $g_t(A)$ коммутируют с $T_k(s)$ при
всех $k$, то, как легко проверить, $g_t(A):D(A_k)\to D(A_k)$ при
всех $k$, а потому и $g_t(A):D(A)\to D(A)$. Отсюда следует, что
$D(A)$ есть существенная область для генератора $G$ (см. [10],
следствие 3.1.7). С другой стороны, $D(A)$ есть существенная
область для оператора $\overline{\psi(A)}$, причем сужения
операторов $\overline{\psi(A)}$ и $G$ на $D(A)$ совпадают с
$\psi(A)$. Поэтому $\overline{\psi(A)}= G$,
 что и завершает доказательство.

Теорема 4  мотивирует окончательный вариант основного определения.

{\bf 5. Определение [6]}. Под {\it значением функции} $\psi$ из
${\cal T}_n$ на наборе $A=(A_1,\ldots ,A_n)$ коммутирующих
операторов из ${\rm Gen}(X)$ будем понимать генератор полугруппы
$g(A)$. Это значение мы далее обозначаем $\psi(A)$.  Возникающее
функциональное исчисление будем называть {\it многомерным
исчислением Бохнера-Филлипса}, или {\it ${\cal T}_n$-исчислением.}

Следующая теорема обобщает на многомерный случай одно утверждение
из [11].

 {\bf 6. Теорема}. {\it Предположим, что полугруппы
$T_j$ сжимающие и удовлетворяют условию
$$
\sum\limits_{j=1}^n\overline{\lim\limits_{t \to
+0}}\|I-T_j(t)\|<2.
$$
Тогда для любой функции $\psi$ из ${\cal T}_n$  оператор $\psi(A)$
является генератором голоморфной полугруппы.}

Доказательство. Не нарушая общности будем считать, что
$c_0=\psi(0)=0$. Положим $b_j=\overline{\lim}_{t \to
+0}\|I-T_j(t)\|$ и выберем $\epsilon > 0$ таким, что $\sum_{j=1}^n
b_j+\epsilon <2$. Найдется такое $\delta >0$, что
$\|I-T_j(t)\|<b_j+\epsilon /n$ при всех $j=1,\ldots ,n; t\in
[0;\delta)$.

Далее, заметим, что при $n>1$
$$
I-T(u)=I-T_1(u_1)+T_1(u_1)(I-T_2(u_2)\ldots T_n(u_n)),
$$
\noindent и значит
$$
\|I-T(u)\|\leq \|I-T_1(u_1)\|+\|(I-T_2(u_2)\ldots T_n(u_n))\|.
$$
\noindent Отсюда по индукции получаем, что
$$
\|I-T(u)\|\leq \sum\limits_{j=1}^n\|I-T_j(u_j)\|,
$$
\noindent а потому при $u\in[0;\delta)^n$ справедливо неравенство
$\|I-T(u)\|\leq \sum_{j=1}^n b_j+\epsilon$. Следовательно, если
$x\in X, \|x\|\leq 1$, то
$$
\|(I-g_t(A))x\|\leq
\int\limits_{\Bbb{R}_+^n}\|I-T(u)\|d\nu_t(u)\|x\|\leq
$$
$$
\leq
\int\limits_{[0;\delta)^n}\|I-T(u)\|d\nu_t(u)+\int\limits_{\Bbb{R}_+^n\setminus
[0;\delta)^n}\|I-T(u)\|d\nu_t(u)\leq
$$
$$
\leq\sum_{j=1}^n b_j+\epsilon+\int\limits_{\Bbb{R}_+^n\setminus
[0;\delta)^n}\|I-T(u)\|d\nu_t(u).
$$
То есть
$$
\|I-g_t(A)\|\leq \sum_{j=1}^n
b_j+\epsilon+\int\limits_{\Bbb{R}_+^n\setminus
[0;\delta)^n}\|I-T(u)\|d\nu_t(u).\eqno(2.1)
$$

 Заметим теперь, что направленность мер $\nu_t$ узко
сходится к мере Дирака $\delta_0$ при $t\to +0$. В самом деле,
преобразование Лапласа ${\cal L}\nu_t(s)=e^{t\psi(s)}$ непрерывно
в точке $s=0$ и ${\cal L}\nu_t(s)\to 1={\cal L}\delta_0(s)$ при
$t\to +0$. Поэтому узкая сходимость вытекает из теоремы
непрерывности для многомерного преобразования Лапласа (см.,
например, [12], гл. IX, \S 5, теорема 3 с). Но так как полугруппы
$T_j$ становятся равномерно непрерывными, ограниченная
 функция $u\mapsto\|I-T(u)\|$ непрерывна на $\Bbb{R}_+^n\setminus
[0;\delta)^n$. Следовательно, переходя в (2.1) к верхнему пределу
при $t\to +0$, получим $\overline{\lim}_{t \to
+0}\|I-g_t(A)\|\leq\sum_{j=1}^n b_j+\epsilon<2$. В силу известного
свойства сильно непрерывных полугрупп  (см., например, [13],
следствие 2.5.7) отсюда следует голоморфность полугруппы $g(A)$,
что и требовалось доказать.

{\bf 7. Следствие}. {\it Пусть пространство $X$ равномерно
выпукло, $T_1$ -- голоморфная полугруппа сжатий в $X$, а операторы
$A_2,\ldots,A_n$ ограничены (если $n>1$). Тогда для любой функции
$\psi$ из ${\cal T}_n$ оператор $\psi(A)$ является генератором
голоморфной полугруппы сжатий.}

Доказательство. Условие теоремы выполнено, поскольку
$\overline{\lim}_{t \to +0}\|(I-T_1(t)\|<2$ (см., например, [13],
следствие 2.5.8) и при $j>1$ справедливы равенства $\lim_{t \to
+0}\|(I-T_j(t)\|=0$.

{\bf 8. Следствие}. {\it Пусть пространство $X$ равномерно
выпукло. Если $T$ -- голоморфная полугруппа сжатий в $X$ с
генератором $A$, то  для любой функции $\psi$ из ${\cal T}_1$
оператор $\psi(A)$ является генератором голоморфной полугруппы
сжатий.}

Для одномерного исчисления справедливо следующее неравенство
моментов.

{\bf 9. Теорема}. {\it Для любого оператора $A\in {\rm Gen}(X)$,
порождающего полугруппу $T$ с оценкой $\|T(u)\|\leq M$, и любой
функции $\psi\in {\cal T}_1$ справедливо неравенство ($x\in D(A),
x\ne 0$)
$$
\|\psi(A)x\|\leq -C_M\psi(-\|Ax\|/\|x\|)\|x\|,
$$
\noindent где} $C_M=(M+1)/(1-e^{-(M+1)/M})$.

Доказательство. Можно считать, что в (1.1) $c_0=c_1=0$. Можно
считать также, что $\|x\|=1$. Формулы (1.2) и (1.1) показывают,
что достаточно доказать неравенство
$$
\|(T(u)-I)x\|\leq C_M \left(1-e^{-\|Ax\|u}\right),\quad u>0.
$$

Обозначим через $t(r)$ функцию, обратную возрастающей функции
$r(t)=t/(1-e^{-t}),\ r(0)=1$. Для фиксированного $r\geq 1$
возможны два случая.

1) $\|Ax\|u\leq t(r)$. Тогда $r(\|Ax\|u)\leq r(t(r))=r$, т. е.
$\|Ax\|u \leq r\cdot\left(1-e^{-\|Ax\|u}\right)$, а потому
$$
\|(T(u)-I)x\|=\|\int\limits_0^1\frac{d}{ds}T(us)xds\|=
$$
$$
\|\int\limits_0^1T(us)Axuds\|\leq M\|Ax\|u\leq
Mr\cdot\left(1-e^{-\|Ax\|u}\right).
$$

2) $\|Ax\|u>t(r)$. Тогда
$$
\|(T(u)-I)x\|\leq
M+1\leq\frac{M+1}{1-e^{-t(r)}}\left(1-e^{-\|Ax\|u}\right).
$$
В любом случае справедливо неравенство
$$
\|(T(u)-I)x\|\leq C(r)\left(1-e^{-\|Ax\|u}\right),
$$
\noindent где $C(r)=M\max\{r;(M+1)/M\left(1-e^{-t(r)}\right)\}$.
Для минимизации $C(r)$ заметим, что функция $t(r)$ возрастает от 0
до $+\infty$ при $1\leq r<+\infty$. Поэтому уравнение
$r=(M+1)/M\left(1-e^{-t(r)}\right)$, т. е. $t(r)=(M+1)/M$, имеет
единственное решение $r_0=r((M+1)/M)$. Если $r<r_0$, то в силу
отмеченной монотонности
$$
\frac{M+1}{M\left(1-e^{-t(r)}\right)}>
\frac{M+1}{M\left(1-e^{-t(r_0)}\right)}=\frac{M+1}{M\left(1-e^{-(M+1)/M)}\right)}=r_0.
$$
Таким образом,
$$
\min\{C(r):r\geq 1\}=Mr_0=(M+1)/\left(1-e^{-(M+1)/M}\right),
$$
\noindent и теорема доказана.

{\bf 10. Следствие}. {\it  Если функция $\psi\in {\cal T}_1$
ограничена на $(-\infty;0)$, то оператор $\psi(A)$ ограничен при
всех $A\in {\rm Gen}(X)$. Обратно, если в некотором банаховом
пространстве $X$ оператор $\psi(A)$ ограничен при всех $A\in {\rm
Gen}(X)$, то  функция $\psi\in {\cal T}_1$ ограничена на
$(-\infty;0)$.}

Доказательство. Первое утверждение очевидно. Для проверки второго
достаточно взять $A=sI$, где $s<0$, и заметить, что тогда
$\psi(A)x=\psi(s)x$.

{\bf 11. Следствие}. {\it Если последовательность функций
$\psi_n\in {\cal T}_1$ сходится к нулю поточечно на $(-\infty;0]$,
то $\psi_n(A)x\to 0$ при всех $x\in D(A)$.}

\bigskip

\begin{center}
ЛИТЕРАТУРА
\end{center}

 1. Kishimoto, A. Subordinate semigroups and order
properties / A. Kishimoto, D. Robinson//J. Austral. Math. Soc.
(Series A) - 1981. - Vol. 31. -- P. 59-76.

2. Berg, C. Generation of generators of holomorphic semigroups /
C. Berg at al. // J. Austral. Math. Soc. (Series A) -- 1993. --
Vol. 55. -- P. 246 - 269.

3. Carasso, A. S. On subordinated holomorphic semigroups /A. S.
Carasso and  T. Kato //  Trans. Am. Math. Soc. --1991. -- Vol.
327. -- P. 867-878.

4.  Миротин, А. Р.  Действие функций класса Шенберга ${\cal T}$ на
конусе диссипативных элементов банаховой алгебры /А. Р. Миротин //
Мат. заметки. -- 1997. -- Т. 61, N 4. -- C. 630 -- 633; English
translation, Math. Notes 61:4 (1997), 524–527.

5. Миротин, А. Р. Функции класса Шенберга ${\cal T}$ действуют в
конусе диссипативных элементов банаховой алгебры, II /А. Р.
Миротин // Мат. заметки. -- 1998. -- Т. 64, N 3. -- C. 423 -- 430; English translation, Math. Notes 64:3 (1998), 364–370.

6.  Миротин, А.Р. Многомерное ${\cal T}$-исчисление от генераторов
$C_0$-полугрупп /А. Р. Миротин // Алгебра и анализ. -- 1999. -- Т.
11, N 2. -- С. 142 -- 170; English translation, St. Petersburg Math. J. 11:2
(2000), 315–335.

7. Пустыльник, Е. И. О функциях позитивного оператора / Е. И.
Пустыльник // Мат. сборник. -- 1982. -- Т. 119, № 1. -- С. 32 --
47.

8. Хилле, Э. Функциональный анализ и полугруппы / Э. Хилле, Р.
Филлипс -- М. : ИЛ, 1962. -- 829 с.

9. Ахиезер, Н. И. Классическая проблема моментов и некоторые
вопросы анализа, связанные с нею / Н. И. Ахиезер -- М.: Физматгиз,
1961. -- 310 с.

10. Браттели, У. Операторные алгебры и квантовая статистическая
механика / У. Браттели,  Д. Робинсон -- М. : Мир, 1982. -- 511 с.

11. Миротин, А.Р. О ${\cal T}$-исчислении  генераторов
$C_0$-полугрупп /А. Р. Миротин // Сибирский мат. ж. -- 1998. --Т.
39, № 3. -- С. 571 -- 582; English translation, Siberian Math. J. 39:3 (1998), 493–503;
“Letter to the editor”, Sibirsk. Mat. Zh. 41:4 (2000), 960; English translation, Siberian
Math. J. 41:4 (2000), 800.

12. Бурбаки, Н. Интегрирование. Меры на локально компактных
пространствах. / Н. Бурбаки -- М. : Наука, 1977. -- 600 с.

13. Pazy, A. Semigroups of Linear Operators and Applications to
Partial Differential Equations/ A. Pazy -- N.Y.: Springer-Verlag,
1983. -- 472 p.

14. Mirotin, A.R. On multidimensional Bochner-Phillips functional calculus, Problemy Fiziki, Matematiki i Tekhniki (Problems of Physics, Mathematics and Technics). -- 2009. --  no. 1. -- p.
60 -- 63 (Russian).

\end{document}